\documentclass[10pt]{article}
\usepackage{graphicx}
\usepackage{amsmath, amsfonts, amssymb, amsthm, enumerate}
\usepackage{amsmath, amsfonts, amssymb, amsthm, enumerate}
\newtheorem{theorem}{Theorem}[section]
\newtheorem{lemma}[theorem]{Lemma}
\newtheorem{proposition}[theorem]{Proposition}
\newtheorem{corollary}[theorem]{Corollary}
\theoremstyle{definition}
\newtheorem{remark}[theorem]{Remark}
\newtheorem{example}[theorem]{Example}

\newtheorem{definition}[theorem]{Definition}

\newcommand{\R}{\mathbb{R}}

\newcommand{\N}{\mathbb{N}}

\renewcommand{\P}{\mathbb{P}}
\newcommand{\T}{\mathcal{T}}

\renewcommand{\d}{\, \mathrm{d}}

\newcommand{\lf}{\left}
\newcommand{\rt}{\right}
\newcommand{\vphi}{\varphi}

\newcommand{\lbl}[1]{\label{#1}}

\newcommand{\eq}[1]{\begin{equation}#1\end{equation}}

\newcommand{\mcal}[1]{\mathcal{#1}}

\makeatletter
\def\moverlay{\mathpalette\mov@rlay}
\def\mov@rlay#1#2{\leavevmode\vtop{%
   \baselineskip\z@skip \lineskiplimit-\maxdimen
   \ialign{\hfil$\m@th#1##$\hfil\cr#2\crcr}}}
\newcommand{\charfusion}[3][\mathord]{
    #1{\ifx#1\mathop\vphantom{#2}\fi
        \mathpalette\mov@rlay{#2\cr#3}
      }
    \ifx#1\mathop\expandafter\displaylimits\fi}
\makeatother

\usepackage{pifont}
\usepackage[numbers,sort&compress]{natbib}

\usepackage[bookmarks,bookmarksnumbered,%
    allbordercolors={0.8 0.8 0.8},%
    linktocpage%
    ]{hyperref}
    \usepackage{bbm}
    
\usepackage{mathtools}

\def\d{\, {\rm d}}

\def \A {\mathcal A}

\def \P {\mathcal P}
\def\T{\mathcal T}

\def\M{\mathcal M}

\newcommand{\mK}{\ensuremath{\mathcal{K}}}

\def\<{\langle}
\def\>{\rangle}

\DeclareMathSymbol{\Z}{\mathbin}{AMSb}{"5A}

\def\d#1{{\rm d}\hbox{$\mskip 0.5mu$}#1}

\title{
{\Large A General Aubry-Mather Theory}}

\makeatletter
\newcommand*\l@chapter[2]{%
  \ifnum \c@tocdepth >\m@ne
    \addpenalty{-\@highpenalty}%
    \vskip 1.0em \@plus\p@
    \setlength\@tempdima{1.5em}%
    \begingroup
      \parindent \z@ \rightskip \@pnumwidth
      \parfillskip -\@pnumwidth
      \leavevmode \bfseries
      \advance\leftskip\@tempdima
      \hskip -\leftskip
      #1\nobreak\hfil \nobreak\hb@xt@\@pnumwidth{\hss #2}\par
      \penalty\@highpenalty
    \endgroup
  \fi}
\renewcommand*\l@section{\@dottedtocline{1}{1.5em}{2.3em}}
\makeatother
\begin{document}

\author{\bf
Nassif  Ghoussoub
\\ \\
{\it Department of Mathematics,  The University of British Columbia}\\
{\it Vancouver, BC, Canada V6T 1Z2}\\
}
\date{}
\maketitle

\begin{abstract}
This paper reproduces the front matter --- preface, overview and table of
contents --- of a monograph by the author, submitted for publication under the title
{\it Skew Linear Entropies and Kantorovich Operators: A General Aubry-Mather
Theory}. The book isolates a class of non-linear operators, which we call
{\it Kantorovich operators}, that are ubiquitous in analysis, probability,
dynamical systems, mathematical economics and finance. 
We develop aspects of their ergodic theory in a way 
that extends classical ones involving Markov operators, free-energy transfers, or the Hopf--Lax--Oleinik semi-group. 
Having no adjoint, the duality between such an operator and measures is
carried instead via a convex functional on {\it pairs} of probability
distributions --- a source and a target --- which we call a {\it skew-linear
entropy}, and which is a general form of optimal mass transport. The extensive overview
reproduced here describes the resulting ergodic theory, in which minimal
measures, a Mather constant, weak KAM solutions and an Aubry set are attached
to an arbitrary Kantorovich operator, extending Aubry--Mather theory well
beyond its origins in Hamiltonian dynamics.
\end{abstract}

\makeatletter
\begingroup
\def\newlabel#1#2{\global\expandafter\def\csname r@#1\endcsname{#2}}

\newlabel{aubry.basic}{{14.2.5}{237}{}{theorem.14.2.5}{}}
\newlabel{ch:Hamiltonian_dynamics}{{22}{359}{Mather Theory in Hamiltonian dynamics}{chapter.22}{}}
\newlabel{ch:risk}{{24}{395}{Controlled and risk-sensitive weak KAM theory}{chapter.24}{}}
\newlabel{chap11_Deterministic.Lagrangian}{{11}{169}{Skew-linear entropies in Lagrangian dynamics}{chapter.11}{}}
\newlabel{chap15_Balayage}{{17}{293}{Weak KAM operators in linear and convex analysis}{chapter.17}{}}
\newlabel{chap16_Complex}{{18}{305}{Weak KAM operators in complex analysis}{chapter.18}{}}
\newlabel{chap1_Linear_Transfers}{{2}{15}{Dual representations of Kantorovich operators}{chapter.2}{}}
\newlabel{chap4_Balayage}{{6}{87}{Balayage and homogeneous Kantorovich operators}{chapter.6}{}}
\newlabel{chap4_convex_capacities}{{3}{39}{Kantorovich operators as functional capacities}{chapter.3}{}}
\newlabel{chap5_Pressures}{{5}{69}{Transport and pressure with cost and control}{chapter.5}{}}
\newlabel{chap5_Stoch_MassTransport}{{12}{189}{Skew-linear entropies in stochastic dynamics}{chapter.12}{}}
\newlabel{continuous_section}{{14}{233}{Mather constants and weak KAM solutions}{chapter.14}{}}
\newlabel{convex_transfers}{{9}{131}{Skew-convex functionals and entropies}{chapter.9}{}}
\newlabel{ergodic_optimiz}{{23}{377}{Mather theory in Ergodic Optimization}{chapter.23}{}}
\newlabel{inequalities_chapter}{{10}{151}{Inequalities between skew-convex entropies}{chapter.10}{}}
\newlabel{multilinear}{{25}{403}{Appendix: Towards a notion of skew-multilinear entropies}{chapter.25}{}}
\newlabel{operations_chapter}{{8}{117}{Operations on skew-linear entropies}{chapter.8}{}}
\newlabel{optimalbalayagetransfer}{{7}{107}{Skew-linear entropies as optimal cost of balayage}{chapter.7}{}}
\newlabel{simplices}{{4}{59}{Kantorovich operators on the simplex}{chapter.4}{}}
\newlabel{skorohod}{{19}{319}{Weak KAM operators and optimal Brownian stopping}{chapter.19}{}}

\endgroup
\makeatother

\section*{Preface}\lbl{preface}

{\it Mathematics is the art of giving the same name to different things.}\hfill Henri Poincar\'e.\vskip 20pt

Many of the non-linear operators encountered in mathematical analysis share several  features. This monograph identifies that common structure and gives it a name: {\it Kantorovich operators}.
These non-linear extensions of \textit{Markov operators} are omnipresent in various fields of mathematical inquiry, albeit under different guises. 

They surface as \textit{convexification operations, Lipschitz extensions, and inf-convolutions} in functional and classical analysis; as \textit{subharmonic and plurisubharmonic envelopes, r\'eduite and optimal stopping operations}  in potential and probability theories; as \textit{non-Markov shifts} in ergodic minimization; 
as \textit{Hopf--Lax--Oleinik semi-groups} in Hamiltonian mechanics; as \textit{risk-sensitive Bellman operators} in operations research and optimal control,
and as \textit{Ruelle transfer operators} in dynamical systems and formal thermodynamics. 
Their name, however, comes from a more specific source: Kantorovich's dual formulation of the Monge optimal mass transport problem.

Kantorovich operators underlie several familiar notions of non-linear analysis.  They induce a special type of {\em functional Choquet capacity}, and inherit from it the monotone continuity needed to make sense of their iterates and limits --- which is why they could have been called {\it convex functional capacities}.   

When a Kantorovich operator's range consists of scalars, it coincides with the \textit{pressure functions} of statistical mechanics and dynamical systems. Integrate a cost function or a control into those pressure functions and the operator version emerges naturally, which is why they might equally have been named \textit{pressure operators}.

Kantorovich operators also have a life in probability theory, where they appear in the study of {\it ``gambling houses''}:  Given a reward function, the maximal expected gain of a gambler --- net of a fee that depends on both their wealth and their chosen distribution of gains --- is exactly a Kantorovich operator. In mathematical finance, they appear as \textit{monetary risk measures}. Kantorovich operators can be seen as a {\it scaled cumulant generating functional} of a \textit{conditional large-deviation principle}. In Markov decision theory, they appear as \textit{risk-sensitive Bellman operators}. 

Our central aim in this book is the ergodic theory of these non-linear operators, and it is here that the classical context of Markov operators is most decisively extended. Mather's theory of minimal measures, refined by Fathi into a  ``weak'' version of Kolmogorov--Arnold--Moser theory, is at heart the ergodic theory of a single non-linear operator --- the Lax--Oleinik operator of Hamiltonian dynamics. A linear operator owes its invariant measures to its adjoint, but a non-linear one has no adjoint. However, Kantorovich operators do possess a dual structure that acts on measures: a rich class of convex functionals defined on pairs of probability distributions --- a source and a target.  We call them {\it skew-linear entropies}.  

This duality extends the classical Legendre duality between pressure functions and general entropies of individual distributions (i.e., convex energy functionals on manifolds of probability densities). The skew-linearity reflects a property of the entropy as a function of one of its two variables.  
It is  ``backward'' (resp., ``forward'') if its Legendre transform as a function of the target (resp., the source) depends linearly on the source (resp., the target). 

 Skew-linear entropies appear classically whenever two probability distributions are linked by a transition probability (Markovian theory), a comparison order (Balayage theory), or a suitably stopped stochastic process (Skorokhod theory). In each of these classical schemes, the associated skew-linear entropy is {\it ``cost-free''}: it takes only the values $0$ and $+\infty$, according to whether the source and target distributions are relatable at all. Cost-free skew-linear entropies are purely qualitative indicators --- yet the Kantorovich operators they generate are far from trivial, though  positively $1$-homogeneous.

General skew-linear entropies integrate costs into these schemes.  Any classical entropy such as the  \textit{logarithmic entropy}, the \textit{Donsker--Varadhan information}, or any  {\it large-deviation rate function} induces a non-trivial, though stationary, skew-linear entropy. So does any point transformation translated by a potential. But it was the optimal mass transport theory of Monge and Kantorovich that led to the introduction of a cost for correlating a pair of distributions, including the above-mentioned classical schemes. However, while the optimal cost of standard mass transportations leads to an entropy that is skew-linear in both directions, forward and backward, its stochastic counterparts, being inherently asymmetric, produce only a one-sided skew-linear entropy and Kantorovich operator.

A key representation theorem shows that most skew-linear entropies are, in fact, the optimal cost of a suitable general mass transport between the two distributions --- one that prices the move of a point mass to a whole distribution, rather than one point to another. Related but broader are the {\it skew-convex entropies}, which include convex functions of skew-linear entropies, %
as well as \textit{logarithmic entropies} when the base measure is considered as a varying source variable.

In our earlier, unpublished preprints, we called skew-linear entropies \emph{linear transfers}, seeing them then as generalizations of optimal mass transport. Our friend and colleague Ivar Ekeland set us straight: these functionals don't describe ``mass transfers'' at all --- they assign a value, or an optimal cost, to a \emph{transfer} of mass. We have since come to see them, more accurately, as extensions of a classical entropy of a single distribution to an ``entropy'' of a pair of distributions. 

Some of this material grew out of collaboration with my former student Malcolm Bowles, and first appeared in his PhD dissertation. This monograph is a considerably expanded and improved version of that joint work. It also draws on joint projects with my colleague Young-Heon Kim, former students Tongseok Lim and Alistair Barton, and postdocs Aaron Palmer and Samer Dweik.

I am grateful to many colleagues and friends for encouraging me to persevere in the challenging task of connecting so many dots across mathematical analysis, and to their universities for the hospitality that made completing it possible: Ivar Ekeland and Eric S\'er\'e at Universit\'e Paris-Dauphine, Changfeng Gui at the University of Macau, Xiaoming Yuan and Dong Li at Hong Kong University, Juncheng Wei at the Chinese University of Hong Kong, and, last but not least, Jingbo Dou at Shaanxi Normal University. I owe particular thanks to Bernard Maurey, who read through several parts of this monograph and supplied a number of illuminating examples, making it more comprehensive, readable, and informative for graduate students and specialists alike. And to Ina Mette, the endlessly patient and friendly acquisitions editor at the American Mathematical Society, for her support throughout.\\

\hfill Nassif Ghoussoub, Vancouver, British Columbia

\bigskip

\section*{Overview}\lbl{introduction}

Mather's theory looks at the invariant measures of a conservative dynamical system and identifies the {\it minimal} ones, those of least time-averaged action. These are the measure-theoretic counterparts of the {\it invariant tori}
of Kolmogorov--Arnold--Moser theory: when the tori of an integrable system break up under perturbation, the minimal
measures are what survives \cite{Evans2}. %
Their least value is a critical constant, the {\it Ma\~n\'e constant};\,   %
they concentrate on a
distinguished invariant set, the {\it Mather set}; and Fathi's weak KAM refinement adds the calibrating potentials, i.e., the fixed points --up to that constant-- of the non-linear 
Lax--Oleinik operator \cite{Mat,Au,Man,Fa}. %
So Mather's theory is, essentially, the {\it ergodic theory of a single non-linear operator}, and the purpose of this book is to present a general framework, where such a theory can be developed further so that it can be applied to a larger class of non-linear operators encountered in mathematical analysis, probability theory, dynamical systems, and formal thermodynamics.

Classical ergodic theory studies a {\it linear} operator $S$ on a space of continuous functions $C(X)$ and draws its invariant measures from the fixed points of the {\it adjoint} $S^*$, which is naturally defined on the dual of $C(X)$. 
These invariant measures, %
together with the leading eigenvalue and eigenfunction of $S$,  govern the long-run dynamics. A {\it non-linear} operator has no adjoint, yet a natural duality with certain actions on measures was present in the work of the above named authors, at least in the case of the Lax--Oleinik operators. This duality was eventually made explicit by Bernard and Buffoni \cite{BB1,BB2} %
 for a class of operators of the form $Tf(x)=\sup\{f(y)-c(x, y); y\in X\}$, which had appeared in the work of Kantorovich \cite{K} 
 as dual objects to the Monge optimal mass transport problem between probability distributions corresponding to a cost $c$ in $ C(X\times X)$.  

This class being still very restrictive to include various operators that appear naturally in probability theory and thermodynamics, we embarked on a decade-long project to find the right encompassing framework. Eventually, a natural notion of non-linear operators emerged for which a general ergodic theory can be developed, with the added bonus that these non-linear extensions of \textit{Markov operators} are omnipresent in various fields of mathematical inquiry, albeit under different guises.  

We call them \emph{Kantorovich operators} in reference to the brief history cited above. They are then in duality with a rich class of convex functionals defined not on single probability distributions but on pairs of them --- a source and a target. 
 We call these {\em skew-linear entropies}. This correspondence extends the classical Legendre duality between pressure functions and general entropies of individual distributions (i.e., convex energy functionals on manifolds of probability densities). The skew-linearity reflects a property of the functional as a function of one of its two variables.  A skew-linear entropy is called {\it forward} or {\it backward} according to whether that variable is the target measure or the source measure. Accordingly, there will be forward or backward Kantorovich operators.

\subsection*{ Kantorovich operators}

To define them, we consider $C(X)$ (resp., $USC(X)$) (resp., $LSC(X)$) to be the space of continuous (resp., the cone of bounded above, proper and upper semi-continuous), (resp., the cone of bounded below, proper and lower semi-continuous) functions on a compact metric space $X$. While a {\it Markov operator} is a positive bounded linear operator $T: C(Y) \to C(X)$  such that $T1=1$,  Kantorovich operators are  defined in the following way:  
 
\begin{definition}\lbl{Kantorovich.op}
\noindent {\it A backward %
Kantorovich operator} is a map $T^-: C(Y) \to USC(X)$ 
that satisfies the following properties: It is  
  \begin{enumerate}
 \item %
   {\it monotone increasing}, i.e., if $g_1, g_2 \in C(Y)$ and $g_1\leq g_2$,  %
  then $T^-g_1\leq T^-g_2$. %
 \item %
 {\it translation invariant on constants}, i.e., if $c\in \R$ and $g\in C(Y)$, then %
$$T^-(g+c)=T^-g +c.$$
 \item  {\it convex}, i.e., 
 if $\lambda \in [0, 1]$, and $g_1, g_2 \in C(Y)$, then 
\begin{equation*}
T^-(\lambda g_1+(1-\lambda)g_2)\leq \lambda T^-g_1+(1-\lambda)T^-g_2.
\end{equation*}
\item {\it lower semi-continuous}, i.e., 
if $g_n \to g$ in $C(Y)$,  %
 then $T^-g\leq \liminf\limits_{n \to \infty}T^-g_n $.
 \end{enumerate}
{\it A forward Kantorovich operator} is a map $T^+: C(X) \to LSC(Y)$ that satisfies 1), 2), with 3) and 4) replaced by 
\begin{enumerate}\setcounter{enumi}{4}
\item  {\it concave}, i.e., 
 if $\lambda \in [0, 1]$, and $f_1, f_2 \in C(X)$, then 
\begin{equation*}
T^+(\lambda f_1+(1-\lambda)f_2)\geq \lambda T^+f_1+(1-\lambda)T^+f_2.
\end{equation*}
\item  \textit{upper semi-continuous} i.e.,  if 
$f_n \to f$ in $C(X)$, 
 then $\limsup_{n \to \infty}T^+f_n \leq T^+f$.
\end{enumerate} 
\end{definition}
Non-linear Kantorovich operators and their iterates already appear -though implicitly- even in classical studies of linear operators. Indeed, if $S: C(X) \to C(X)$ is a Markov operator, then it is clearly a Kantorovich operator, but so is
any operator of the form 
\eq{T f= Sf-A,}
and its iterates, where 
$A$ is a continuous ``observable" considered in ergodic optimization models (See Chapter \ref{ergodic_optimiz}). 

The free-energy operator \eq{Tf=\log S(e^f)} associated to $S$ and its iterates, considered in formal thermodynamics, are also backward Kantorovich operators, often used to analyze the spectral properties of $S$. 

Even more basic are the maps    \eq{T^- g:=g\vee Sg  \quad \hbox{(resp.,\,\, $T^+f:=f\wedge Sf$}),}
 which are also non-linear backward (resp., forward) Kantorovich operators. 
Iterating $T^-$ and $T^+$ often lead to the  {\it idempotent} backward  (resp., forward) Kantorovich operator 
 \eq{T^-_\infty  g=\lim_n\uparrow T^-_ng \quad \hbox{(resp., 
 $T^+_\infty f=\lim_n\downarrow T^+_nf$),}
 }
  which are none other than the 
  {\it r\'eduite}  operators
 \eq{T^-_\infty g (x)=\hat g (x):=\inf \left\{h(x); -h \in \A, h\geq g \,\, {\rm on}\,\, X\right\}
 }
 and 
 \eq{T^+_\infty f(x)=\check{f}(x):=\sup \left\{h(x); h \in \A, h\leq f \,\, {\rm on}\,\, X\right\},}
 corresponding  to the cone  $
\A=\{f\in C(X); f\leq Sf\}
$
of $S$-subharmonic functions, yielding %
 the least $S$-superharmonic function above $g$ (resp., the greatest $S$-subharmonic function below $f$). 
Note that in both examples  
\eq{
  T^-\circ T^-_\infty  g=T^-_\infty  g\quad {\rm and}\quad T^+\circ T^+_\infty g= T^+_\infty g,
  } which means that the ranges of the operators $T^-_\infty$ and $T^+_\infty$ are the sets of fixed points of $T^-$ and $T^+$ respectively. 

Actually, one need not start with a Markov operator to construct envelopes. Indeed, if $X$ is a strictly convex open domain in $\R^d$, then the {\it concavification operator} (Chapter \ref{chap15_Balayage})
\eq{
T g(x)=\sup \{\frac{1}{2}(g(x_1)+g(x_2)); x_1, x_2 \in X, x=\frac{1}{2}(x_1+x_2)\}
}
  is a backward Kantorovich operator 
  whose iterates induce an idempotent Kantorovich operator that associates to a function, its upper semi-continuous concave envelope, equivalently, 
 \eq{
 T_{\infty}g(x)=\sup\left\{ \mathbb{E}[g (F_n)]; n\in \N, (F_k)^n_{k=0} \hbox{ dyadic martingale with $F_0 = x$} \right\}.
 }
 In the theory of several complex variables (Chapter \ref{chap16_Complex}),  the operator
 \eq{
T g(x):= \sup_{v\in \R^n} \bigg\{ \int ^{2\pi}_0 g(x +
 e^{i\theta}v) {d\theta \over 2\pi};\,  x + \bar \Delta v\subset O \bigg\}, 
 }
where $O$ is a bounded strictly pseudoconvex domain in ${\mathbb C}^d$, $\Delta = \{ z \in {\mathbb C}, \vert z\vert < 1\}$ is the open unit disc in ${\mathbb C}$, is a Kantorovich operator  and its iterates lead to the idempotent backward Kantorovich operator giving the  {\it pluri-superharmonic envelope} of functions on $\overline O$ %
 \eq{
T_{\infty}g(x)=\sup \left\{ \int^{2\pi}_0 g (P(e^{i\theta})) {d\theta
\over 2\pi};\ \hbox{\rm $P$ polynomial, $P(\bar \Delta)\subset O$,  
 $P(0) = x$}\right\}.
 }
If $O$ is a convex bounded domain in ${\mathbb R}^n$, the operator  
 \eq{
T g(x) =\sup_{r \geq 0} \bigg\{ \int_{B} g (x + r y) \,dm(y);\,  x + r \overline{B} \subset O \bigg\},
}
where $B$ is the open unit ball in $\R^n$ centered at $0$, and $m$ is normalized Lebesgue measure on $\R^n$, is a Kantorovich operator and its iterates lead to the {\it superharmonic envelope via optimal stopping} %
 \eq{T_{\infty} g(x)
		:= \sup\Big\{{\mathbb E}^{x}\Big[g(B_\tau)\Big]; \, \tau \ge 0 \,\, \hbox{ stopping\, time in ${\mathcal S}(O)$} %
		\Big\}.
}
The expectation $\mathbb{E}^{x}$  refers to Brownian motions $(B_t)_t$ starting at $x$ and ${\mathcal S}(O)$ is the class of (possibly randomized) Brownian stopping times $\tau$ prior to the exit time  $\tau_O:=\inf\{t;\ B_t\not\in O\}$ (Chapter \ref{skorohod}). 

If $X$ is a Banach space, then regularizing operators such as inf-convolution 
\eq{
Tf(x)=\inf\{f(y)+\|x-y\|^2; y\in X\},
}
 (resp., Lipschitz extension of functions on a subset $K \subset X$, namely
 \eq{
 Tg (x)=\sup\{g(y)-\|x-y\|; y\in K\}),
 }
 are typical  forward (resp., backward and idempotent) Kantorovich operators.
 
 More generally, if $c:X\times Y\to \R$ is a bounded below lower semi-continuous cost function, then  the operators defined for any $f\in C(X)$ (resp., $g\in C(Y)$) by 
\begin{equation}
T ^+_cf(y)=\inf_{x\in X} \{f(x)+c(x, y)\} \quad {\rm and} \quad T ^-_cg(x)=\sup_{y\in Y} \{g(y)-c(x, y)\},
\end{equation}
which appear in the Monge-Kantorovich duality of optimal mass transport associated to the cost $c$, are also forward (resp., backward) Kantorovich operators (Chapter \ref{chap5_Pressures}).  

One can also consider {\it balayage operations with cost} pointing  to elements of a``non-linear potential theory". For example,  if $O$ is a convex bounded domain in $\R^d$, then the operator 
$Tg(x)=u_{g, x},$
 where $u_{g, x}$ is the unique minimiser of the variational problem
\eq{\lbl{gen.cap}
	\inf\Big\{\int_O \big|\nabla u \big|^2dy;\ u\geq g- c(x,\cdot),\,  u\in H^1(O)\Big\},   
}
is also a backward  Kantorovich operator.

{\it Pressure functions} (i.e., essentially Legendre conjugates of energy functionals or general entropies) are Kantorovich operators whose range consists of constant functions. There are two mechanisms that turn pressure functions into full-fledged Kantorovich operators: one couples the base to the target through a cost, and another aggregates over a control. For example, starting with  the logarithmic entropy pressure $Pg= \log \int_{Y}e^{g}d\nu
$ as ``a seed", we can subject it  to a cost $c$ to obtain the Sinkhorn operator 
\eq{\lbl{entro.0}
T_cg(x)= \epsilon \log \int_{Y}e^{\frac{g(y)-c(x,y)}{\epsilon}}d\nu(y),
}
which figures in the {\it entropic regularization} of optimal mass transports (Chapter \ref{operations_chapter}).

Another mechanism is to consider suitably controlled entropies such as in the foundational \emph{risk-sensitive Bellman operator} in the theory of \emph{dynamic programming},
\eq{\lbl{risk.bellman0}
T_\gamma g(x)=\max_{a\in\mcal A}\Big\{\,r(x,a)+\tfrac1\gamma\log\int_X e^{\gamma g(y)}\,P_a(x,dy)\Big\},
\qquad g\in C(X),
}
where $\mcal A$ is a compact set of controls, $r\in C(X\times\mcal A)$ is a reward function, %
and $P_a(x,\cdot)\in\mcal P(X)$ are transition kernels  depending %
continuously on $(x,a)$ (Chapter \ref{chap5_Pressures}). This Kantorovich operator  leads to many other known ones according to when we turn off \emph{the risk} ($\gamma\to0$, i.e., average-cost Bellman), the temperature ($\gamma\to\infty$, i.e.,  robust control), or the control (i.e., $\A$ is a singleton). In the latter case, and  when $\gamma=1$, we obtain the celebrated free energy Ruelle operator. The most degenerate case, when $\gamma\to0$ and $\mcal A$ is a singleton are the operators appearing in ergodic minimization (Chapters \ref{ergodic_optimiz}
and \ref{ch:risk}). 

The logarithmic entropy seed pressure is only one choice. The control construction %
applies to other pressure seeds, and each produces a Bellman-type operator in which a controller optimizes that pressure. One can then obtain a rich class of Kantorovich operators besides the Bellman operator, such as the {\em controlled thermodynamic formalism}, the {\em controlled topological pressure}, and the {\em controlled Fisher--Donsker--Varadhan information}. 

The framework of this monograph originated in Hamiltonian dynamics.  If $L: TM \to \R$ is a given \textit{Tonelli Lagrangian} on a smooth compact Riemannian manifold $M$ without boundary, %
 and $H$ is the associated Hamiltonian, then the ``Hopf-Lax-Oleinik" operator 
\begin{equation}\label{value.20}
T^-g(x) :=\sup\Big\{g(\gamma (1))-\int_0^1L(\gamma (s), {\dot \gamma}(s))\, ds; \gamma \in C^1([0, 1], M);   \gamma(0)=x\Big\},
\end{equation}
which associates to a state $g$ at time $0$, the state at time $1$ of the  viscosity solution for the associated backward  Hamilton-Jacobi equation,  
 \eq{\label{HJ.0} 
\left\{ \begin{array}{lll}
\partial_tV-H(x, \nabla_xV)&=&0 \,\, \text{ on }\,\,  (0, 1)\times M\\
\hfill V(0, x)&=&g(x),
\end{array}  \right.
}
is a backward Kantorovich operator.
Similarly, 
\eq{\lbl{Tplus}
T^+ f(y) := \inf\{ f(\gamma(0)) + \int_{0}^{1}L(\gamma(s), \dot{\gamma}(s))d s\,;\, \gamma \in C^1([0,1]; M), \gamma(1) = y\},
}
is a forward Kantorovich operator, which associates to an initial state $f$, the state  at time $1$ of the viscosity solution for the corresponding forward Hamilton-Jacobi equation (Chapter \ref{chap11_Deterministic.Lagrangian}). Iterates of these operators lead to fixed points -up to an important constant $c$- that are viscosity solutions of the underlying stationary Hamilton-Jacobi equation (Chapter \ref{ch:Hamiltonian_dynamics}).  This important example is at the basis of the work --in the deterministic case-- of  
Mather, 
Aubry 
Ma\~n\'e,  
Fathi, 
 and others 
on relaxed versions of the celebrated Kolmogorov--Arnold--Moser theorem in Hamiltonian dynamics.

\subsection*{Skew-linear entropies of pairs of probability distributions}

 \textit{Duality} is omnipresent in all the examples mentioned above. Indeed, the cone of Kantorovich operators is in natural duality with a rich class of functionals acting on pairs of probability measures, that we call the {\it skew-linear entropies} of a pair of probability distributions. The duality is illustrated by the following crucial representation.

  \begin{theorem}\lbl{repr.op} Let $T$ be a map from $C(Y)$ to the space $B^b(X)$ of Borel functions on $X$ that are bounded above. The following assertions are then equivalent:
  \begin{enumerate}
  \item $T$ is a backward Kantorovich operator from $C(Y)$ to $USC(X)$. 

\item There exists a convex weak$^*$-lower semi-continuous function $\T$ 
on $\P(X) \times \P(Y)$ such that for all $\mu \in \P(X)$ and $g \in C(Y)$,
\eq{
\T_\mu^*(g) = \int_{X}T g\,  d\mu. 
}
Here, $\T_\mu^*$ is the Legendre transform on $C(Y)$ of the convex functional $\T_\mu: \nu \to {\mathcal T} (\mu, \nu)$ constrained to be $+\infty$ outside $\P(Y)$. 

\item There exists a proper lower semi-continuous function $c: X\times {\mathcal P}(Y)\to  \R \cup\{+\infty\}$ such that for all $x\in X$, the functional $\sigma \mapsto c(x, \sigma)$ is convex, and for any $g\in C(Y)$,
\eq{\lbl{4prop.intro}
Tg(x) := \sup\{ \int_{Y}gd\sigma - c(x,\sigma)\,;\, \sigma \in \P(Y) \}.  
}
\end{enumerate}
\end{theorem}
We shall then say that {\it $\T$ is a backward skew-linear entropy}. 

There is a similar statement for a forward Kantorovich operator $T^+$, when one varies the source measure so that 
\eq{
\T_\nu^*(f) = -\int_{Y}T^+ (-f) d\nu \quad \hbox{for $f\in C(X)$,} 
}
where $\T_\nu: \mu \to {\mathcal T} (\mu, \nu)$. $\T$ will then be a {\it forward skew-linear entropy}. 

The duality  
is then illustrated by the formulae: For all $\mu \in \mcal{P}(X), \nu \in \mcal{P}(Y)$,
\eq{
\T(\mu,\nu) = \sup_{g \in C(Y)}\{\int_{Y}gd\nu - \int_{X}T^-g d\mu\}
}
resp.,
\eq{\T(\mu,\nu) = \sup_{f \in C(X)}\{\int_{Y}T^+fd\nu - \int_{X}f d\mu\}.
}

This duality has a probabilistic context that recurs throughout the book. A convex energy $I$ on ${\cal P}(Y)$ is a {\it large-deviation rate function}, and its pressure $I^*$ is the {\it scaled cumulant generating functional} that Varadhan's lemma attaches to it. In  Chapter~\ref{chap1_Linear_Transfers}, it is shown that the general, state-dependent duality 
is the {\it conditional} form of this correspondence: The backward Kantorovich operator $T$ is the generating functional of a family of laws conditioned on a base point $x$, i.e.,  they are transition kernels $\mu_n(x,\cdot)\in\mcal{P}(\mcal{P}(Y))$ and $\T$ is the conditional rate function, in such a way that for $g\in C(Y)$,
\eq{\lbl{cond.T0}
Tg(x):=\lim_{n\to\infty}\tfrac1n\log\int_{\mcal{P}(Y)}e^{\,n\int_Y g\,d\sigma}\,\mu_n(x,d\sigma)=\sup_{\sigma\in\mcal{P}(Y)}\Big\{\int_Y g\,d\sigma-\T(x,\sigma)\Big\}.
}
The entropic  (Sinkhorn), thermodynamic (free-energy), and small-noise operators (Lax--Oleinik) are essentially the conditional principles of Sanov, of Kifer, and of Freidlin--Wentzell.

Besides illustrating the duality, Theorem \ref{repr.op} is crucial for the study of the ergodic properties of Kantorovich operators. First, it allows for the extension of $T$ to become a {\it Choquet functional capacity} defined  on the class $F_+(Y)$  of all non-negative functions valued in $\R\cup\{+\infty\}$ (Chapter \ref{chap4_convex_capacities}).
The extension maps $USC(Y)$ to $USC(X)$ hence allowing for the iteration of such an operator (if $Y=X$). Secondly, the capacity properties of the extension justifies 
monotonic  pointwise limits such as:
\eq{
T(\lim_n\downarrow T^n(g +nc))+c=\lim_n \downarrow T^{n+1}g +(n+1)c.
} 
Mathematical and statistical analysis feature many examples of backward skew-linear entropies. 
For example, the skew-linear entropy associated to 
the Markov operator $S$, is
\begin{equation}
{\mathcal T}(\mu, \nu)=\left\{ \begin{array}{llll}
0 \quad &\hbox{if $\nu = S^*(\mu)$}\\
+\infty \quad &\hbox{\rm otherwise,}
\end{array} \right.
\end{equation}
where $S^*:{\mathcal M}(X) \to {\mathcal M}(X)$ is the adjoint operator. The backward and forward skew-linear entropy associated to the Kantorovich operators $T^-_\infty$ and $T^+_\infty$ considered above  is given by 
\begin{equation}
{\mathcal T}(\mu, \nu)=\left\{ \begin{array}{llll}
0 \quad &\hbox{if  $\mu \prec_{\mathcal A} \nu$ }\\
+\infty \quad &\hbox{\rm otherwise,}
\end{array} \right.
\end{equation}
where ${\mathcal A}$ is the cone of $S$-subharmonic functions and $ \prec_{\mathcal A}$ is the ${\mathcal A}$-\textit{balayage} partial order 
between probability measures $\mu$, $\nu$ in $\P(X)$, that is
 \eq{\lbl{true.0}
 \mu \prec_{\mathcal A} \nu \quad \hbox{ if and only if \quad $\int_X\vphi \, d\mu \leq \int_X\vphi \, d\nu$ for all $\vphi$ in  ${\mathcal A}$.}
 }
There are many other balayage schemes in analysis, 
such as the convex order, the subharmonic order and the plurisubharmonic orders between measures (Chapter \ref{chap4_Balayage}). 
These classical schemes do not involve any cost for their implementation. Their corresponding skew-linear entropy will be $0$ or $+\infty$ according to whether they are related or not. Their associated Kantorovich operators are then {\it positively $1$-homogeneous}, that is they satisfy
\eq{
T(\lambda f)=\lambda Tf \hbox{ for any continuous $f$ and any $\lambda \in \R^+$}.
}

A simple example of a non-zero-cost backward skew-linear entropy is induced by the {\it classical entropy}. By that, we mean any lower semi-continuous convex energy $I$ on  ${\cal P}(X)$. It then induces  
a {\it stationary} backward skew-linear entropy,
\eq{
\T(\mu, \nu)=I(\nu) \quad \hbox{for all $\mu \in {\cal P}(X)$,}
}
whose associated Kantorovich operator is given by the Legendre transform $f\to I^*(f)$, i.e. {\it a pressure function}.

The most concrete examples of skew-linear entropies can be obtained by restricting the underlying probability spaces to $n$-simplices, i.e., when 
$$
   P(X) 
 = \{ \alpha = (\alpha_1, \alpha_2, ...\alpha_n) : 
       \alpha_i \ge 0, \, \alpha_1 + \alpha_2+...\alpha_n = 1 \}.
$$
$C(X)$ can then be identified with $\R^n$ and the duality between $C(X)$ and $\M(X)$ is given by
$
 x \cdot \alpha = x_1 \alpha_1 + x_2 \alpha_2 +...+x_n\alpha_n.
$
Typical examples of backward skew-linear entropies on $\P(X)\times \P(X)$ are given by 
\eq{
   \T(\beta, \alpha) 
 = \sum_{i=1}^n \beta_i G_i(\alpha_i / \beta_i),
}
with corresponding Kantorovich operator defined for each $x=(x_i)_{i=1}^n\in \R^n$ by 
\eq{
T^-((x_i)_{i=1}^n) = \bigl( x^* + g^*_i(x_i - x^*) \bigr)_{i=1}^n,
}
where $x^*:=\max\limits_{i=1,...n} x_i$, each $G_i$ is an appropriate convex function on $[0, +\infty)$, and $g_i^*$ is the Legendre transform of the function $g_i$ that is equal to $G_i$ on $[0, 1]$ and to~$+\infty$ elsewhere (Chapter \ref{simplices}).
  
It is however optimal mass transport theory that gave the idea of introducing a cost for all schemes relating two probability measures, including the classical ones mentioned above (Chapter \ref{chap5_Pressures}). 
Indeed, 
skew-linear entropies include --and originate from-- the value function in deterministic optimal mass transport: If  $c(x, y)$ 
 is a lower semi-continuous cost function $c(x, y)$ on a product space $X\times Y$, then 
 the following functional on ${\mathcal P}(X)\times {\mathcal P}(Y)$ is  both a backward and a forward skew-linear entropy:
\eq{\lbl{MT.0}
{\mathcal T}_c(\mu, \nu):=\inf\big\{\int_{X\times Y} c(x, y) \, d\pi; \pi\in \mK(\mu,\nu)\big\},
}
where %
$\mK(\mu,\nu)$ is the set of probability measures $\pi$ on $X\times Y$ whose marginal on $X$ (resp. on $Y$) is $\mu$ (resp., $\nu$) {\it (i.e., the transport plans)}. 

On the other hand, the \textit{stochastic counterparts} of optimal mass transport provide one-sided examples of skew-linear entropies, which led us to distinguish between the forward and backward notions (Chapter \ref{chap5_Stoch_MassTransport}). Indeed,  optimal transportation problems for stochastic processes with controlled dynamics and fixed end time, where the transport cost is given by a general Lagrangian $L$ on $\R^+\times \R^d\times \R^d$ are intrinsically asymmetric. The functionals 
\begin{align}\label{Primal with fixed end time.intro}
\T_L (\mu,\nu) :=\inf_{\beta} \Big\{\mathbb{E}\Big[\int_0^1 L\big(t,X_t,\beta_t\big)dt\Big];\ dX_t=\beta_t\, dt+dW_t,\ X_0\sim \mu,\ X_1\sim \nu\Big\},
\end{align}
are backward skew-linear entropies but not forward, and their corresponding 
  Kantorovich operator is given by $T^-g:=J_g(0, x)$, 
  where $J_g$ is the initial state of the backward second order Hamilton-Jacobi equation 
\begin{equation} \label{PDE with fixed end time.0}
\begin{cases}
\partial_t J(t,x)  + \frac{1}{2} \Delta J(t,x) + H\big(x,\nabla J(t,x)\big)&= 0 \  \mbox{in\ }\, (0,1) \times \mathbb{R}^d,\\
\hfill J(1,x)&= g(x) \  \mbox{on\ }\,\mathbb{R}^d.
\end{cases}
\end{equation}
Optimal Brownian stoppings with Monge-type cost
 \begin{align}
	\mathcal{B}_c(\mu, \nu) = \begin{cases} \inf%
	\Big\{\mathbb{E} \big[ c(B_0, B_\tau)\big]; \ B_0 \sim \mu,\, B_\tau \sim \nu; \tau \in {\mathcal S}(O)\Big\}\\
	+ \infty \quad \hbox{if no such stopping time exists,} 
	\end{cases}
\end{align}
and optimal Brownian stoppings with Lagrangian -type cost, that is 
	\begin{align} \label{eqn:Skorokhod_cost.0}
		{\mathcal P}_L(\mu,\nu) := \begin{cases} \inf
		\Big\{\mathbb{E}\Big[ \int_0^\tau L(t,B_t)dt\Big];\ B_0 \sim \mu,\, B_\tau \sim \nu; \tau \in {\mathcal S}(O)\Big\}\\
	+ \infty \quad \hbox{if no such stopping time exists,} 
	\end{cases}
\end{align}
provide other examples of one-sided skew-linear entropies and Kantorovich operators.

It turns out that skew-linear entropies can always be represented as the values of an optimal mass transport provided we replace the cost $c(x, y)$ on $X\times Y$ above by a more general cost $c(x,\sigma)$ on $X\times \P(Y)$, and when the choice $(x, \sigma)$ follows a rule given by a balayage cone $\A$ (Chapter \ref{optimalbalayagetransfer}). Indeed,  a standard backward skew-linear entropy $\T$ can always be represented as:
\begin{equation}\label{optimalbalayage.000.intro}
\T(\mu,\nu) = \begin{cases}
\inf\{\int_{X} c(x, \pi_x)d\mu(x)\,;\, \pi \in \mathcal{K}_{\A}(\mu,\nu)\} & \text{if }\mu \prec\prec_{\A} \nu,\\
+\infty & \text{otherwise,}
\end{cases}
\end{equation}
where  $\A$ is an admissible balayage cone on the disjoint union  $X\sqcup Y$, and $(\pi_x)_x$ is a disintegration of $\pi$ with respect to $\mu$ such that $\delta_x\prec\prec_\A \pi_x$ for $\mu$-almost all $x\in X$. This makes a clear connection with the theory of weak mass transports introduced by Gozlan et al.\ \cite{Go4}. %

A natural extension of skew-linear entropies is introduced: the \textit{skew-convex entropies} (Chapter \ref{convex_transfers}). These are functionals $\T(\mu,\nu)$ on $\P(X)\times \P(Y)$ which are suprema of skew-linear functionals $(\T_i)_{i\in I}$ in such a way that 
\eq{ \hbox{$
\T_\mu^*(g)= \inf\limits_{i \in I}\int_{X}T^-_{i}g d\mu$\,\,\, for $\mu \in \P(X)$ and $g \in C(Y)$,
}
} 
 where $(T_i)_{i\in I}$ are (auxiliary) operators associated with $(\T_i)_{i\in I}$.
Note that this implies that the map $\mu \to -\T^*_\mu$ is convex on $\P(X)$ as opposed to being linear. %

An important example of a skew-convex, but not skew-linear entropy is 
the relative logarithmic entropy ${\cal H}(\mu, \nu)$ defined as
\eq{
\hbox{${\cal H}(\mu, \nu):=\int_X\log (\frac {d\nu}{d\mu})\, d\nu$ if $\nu<<\mu$ and $+\infty$ otherwise.}
}
Note that
 \eq{
  {\cal H}_\mu^*(f)=
  \log\int_Xe^f\, d\mu
  }
and $\mu \to  {\cal H}_\mu^*(f)$ is concave and not linear.  On the other hand, convex functions of skew-linear entropies are skew-convex entropies. 

The duality between skew-linear entropies and Kantorovich operators lead to dual statements for various (generalizations of) standard ``transport-entropy'' type inequalities such as the one stating that for a  certain reference probability $\mu$ and a certain cost function $c(x, y)$, one has 
\eq{\T_c(\sigma,\mu) \leq \mcal{H}(\mu,\sigma) \quad \hbox{for all probability $\sigma$}, 
}
if and only if 
\eq{ \int_Xe^{-T_c^-g} d\mu \leq e^{-\int_Xg d\mu} \quad \hbox{for every $g\in C(X)$,}   }
where $\T_c(\mu,\nu)$ is an optimal transport associated to the cost $c$, and $T^-_c$ is the corresponding Kantorovich operator. A more general duality also holds for {\it Maurey-type inequalities}, where in this case, we have two reference measures $\mu, \nu$ in such a way that
\eq{
{\mathcal F}(\sigma_1, \sigma_2) \leq \lambda_1 {\mathcal H}\star {\mathcal T}_1 (\mu,\sigma_1)+\lambda_2 {\mathcal H}\star {\mathcal T}_2( \nu, \sigma_2) \quad \hbox{for all $\sigma_1\in {\mathcal P}(Y_1), \sigma_2\in {\mathcal P}(Y_2)$}
}
if and only if 
\eq{
\sup_{k }\lf(\int_{X_1} e^{T_1^-(- \frac{1}{\lambda_1}F_k^-g)}\, d\mu\rt)^{\lambda_1} \lf(\int_{X_2}e^{T_2^-(\frac{1}{\lambda_2}g)}\, d\nu\rt)^{\lambda_2} \leq 1 \quad \hbox{for all $g\in C(Y_2)$,}
 }
where $\lambda_1, \lambda_2 \in \R^+$, $\mathcal F$ is a skew-convex functional with auxiliary operators $(F_k)_k$, while $\T_1, \T_2$ are backward skew-linear entropies with Kantorovich operators $T_1^-$ and $T_2^-$ (Chapter \ref{inequalities_chapter}).

\subsection*{General weak KAM theories}

One can associate to any backward Kantorovich operator $T$  the following --possibly infinite-- {\it Mather constant}, 
 \eq{\lbl{Mane}
c(T) :=  \inf_{\mu \in \P(X)}\sup_{h\in C(X)}\Big\{\int_X(h-T h)\, d\mu\Big\}.  
}
Dually, 
\eq{\lbl{stc10}
c(T)= \inf_{\mu \in \mcal{P}(X)}\T(\mu,\mu), 
}
where $\T:\mcal{P}(X)\times\mcal{P}(X)\to \R\cup\{+\infty\}$ is the skew-linear entropy associated to $T$. The {\it minimal measures} are those where this infimum is attained. 

The ergodic quantities are then read
off $\T$ rather than off a missing adjoint for $T$:
\begin{itemize}
\item the invariant measures become the {\it minimal} measures.
\item the leading eigenvalue becomes the {\it Mather constant} $c(T)$. %
\item the eigenfunction becomes a {\it weak KAM solution}, a fixed point $Th+c(T)=h$;
\item and occasionally, the recurrent set, carrier of the minimal measures, becomes the {\it Aubry set}.
\end{itemize}
All of this can be read through the conditional large-deviation principle described above, in which $T$ is the
scaled cumulant generating functional and $\T$ the conditional rate function. The Mather constant $c(T)$ is then a
{\it long-time large-deviation rate}; a weak KAM solution $\bar u$ is a {\it tilt}, the exponential change of
measure $e^{\bar u}$ that reweights the conditioned laws onto the optimal behaviour; and a minimal measure is the
{\it typical} one, the law on which the empirical distribution concentrates by the law of large numbers, every
other outcome being exponentially rare (Chapter~\ref{continuous_section}).

When $T$ is itself a Markov operator these reduce to the familiar objects. The Mather constant $c(T)=0$, the minimal measures $\T(\mu,\mu)=0$ are precisely the
invariant measures, and the weak KAM solutions are the invariant functions $Th=h$.     The theory below is what survives once linearity is dropped.

Our terminology is in reference to the above mentioned works  on Hamiltonian dynamics.   
But our main premise is that there are ``several weak KAM theories"  out there  to which this general theory is applicable. For example, in the stochastic counterpart of  Hamiltonian dynamics, in Markov decision theory, in Ruelle's thermodynamic formalism, in mathematical finance, and in the theory of ergodic minimization for expanding dynamical systems.

By a ``weak KAM theory", we mean the process of associating to a backward Kantorovich operator $T$ %
 with a finite Mather  constant $c(T)$, hence to the corresponding skew-linear entropy $\T$, a {\it weak KAM operator} $T_\infty$ %
with its skew-linear entropy $\T_\infty$ ({\it a Peierls barrier}) in such a way that $T_\infty: C(X) \to USC(X)$  is a Kantorovich operator satisfying the following: 
 
 \begin{enumerate}
\item $T_\infty$ is idempotent.

\item $TT_\infty =T_\infty T$.

\item $T_\infty$ maps $C(X)$ to the backward weak KAM solutions for $T$, i.e., for any $g\in C(X)$, 
\begin{equation}
TT_\infty g +c(T)=T_\infty g.
\end{equation}
\item The  associated {\it measure-level Aubry set}, i.e,   
\eq{
\mcal{N}(\T_\infty)=\{\sigma\in\mcal{P}(X);\ \T_\infty(\sigma,\sigma)=0\}, 
}
contains the set  $\mcal M (\T)$ of all minimal measures. This can be seen as the measure-level counterpart to instances when the point-set level Aubry set contains the support of all minimal measures (Proposition \ref{aubry.basic}). 

\item  Every weak KAM solution $\bar u$ is
recovered --via a Peierls barrier $\T_\infty$--  from its values on the corresponding $\mcal{N}(\T_\infty)$ by
\eq{\lbl{aubry.repr.intro}
\int_X\bar u\,d\mu=\sup_{\sigma\in\mcal{N}(\T_\infty)}\Big\{\int_X\bar u\,d\sigma-\T_\infty(\mu,\sigma)\Big\},
\qquad\mu\in\mcal{P}(X).
}
In particular, $\mcal{N}(\T_\infty)$ is a {\it uniqueness set}, i.e.,  two weak KAM solutions with
$\int_X\bar u_1\,d\sigma=\int_X\bar u_2\,d\sigma$ for all $\sigma\in\mcal{N}(\T_\infty)$, are then equal everywhere.

\item  If $\T$ is, in addition, a forward skew-linear entropy, then $\T_\infty$ is also a forward skew-linear entropy with an associated  forward weak KAM operator $T^+_\infty: C(X) \to LSC(X)$ with similar corresponding properties. 
 
 Moreover, there exist a backward (resp., forward) weak KAM operator  $S^-: C(X)\to USC (X)$ (resp., $S^+: C(X)\to LSC(X)$) %
 such that 

 \eq{
S^-=T_\infty^-S^+, \quad S^+=T_\infty^+S^-,
} 
\eq{
T^-S^-+c=S^-, \quad T^+S^+-c=S^+, 
}
and  for any $f\in C(X)$, %
\eq{
\int_X S^+fd\mu=\int_X S^-f d\mu \hbox{\,\,\,  for every  $\mu\in \cal N(\T_\infty)$} %
}
 and for any $\mu, \nu$ in $\P(X)$, we have
\eq{
\T_\infty (\mu, \nu)=\sup\{\int_XS^+f\, d\nu-\int_XS^-f\, d\mu; f\in C(X)\}.
}
\end{enumerate}
These statements are reminiscent of the results of Mather, Fathi, and Bernard-Buffoni when $T$ is the Lax-Oleinik Kantorovich operators.  We shall present  many other situations where such a theory can be applied.

Throughout this monograph, we shall focus  on probability measures on compact spaces, even though the right settings for most applications and examples are complete metric spaces, Riemannian manifolds, or at least $\R^n$. This will allow us to avoid the usual functional analytic complications, and concentrate on the {\it algebraic} aspects of the theory. The simple compact case will at least point to results that can be expected to hold and be proved --albeit with additional analysis and suitable hypotheses -- in more general situations. In the case of complete metric spaces (such as $\R^n$), which is the setting for many examples, one needs to consider more suitable dualities such as  between the space of all bounded and Lipschitz functions and the space of Radon measures with finite first moment. 

Finally, in the last chapter \ref{multilinear}, we explore two different approaches to extend the above theory and define {\it skew-multilinear entropies} on products $ \mcal{P}(X_1)\times \ldots \times \mcal{P}(X_n)$, where $X_j, j = 1,\ldots, n$ are compact spaces. The main example being the multimarginal optimal transport problems with cost function $c : X_1 \times\ldots\times X_n \to \R \cup \{+\infty\}$, 
\eq{
\T_{c}(\mu_1,\ldots,\mu_n) := \inf_{\pi \in \Gamma(\mu_1,\ldots,\mu_n)}\int c\d\pi, 
}
where $\Gamma(\mu_1,\ldots,\mu_n)$ is the set of probabilities on 
$\prod_{j = 1}^{n}X_j$ with marginals on $X_1, X_2,..., X_n$ are  $\mu_1, \mu_2,..., \mu_n$ respectively.

\clearpage

\makeatletter
\renewcommand\contentsline[4]{\csname l@#1\endcsname{#2}{#3}}
\section*{\contentsname}
\@mkboth{\MakeUppercase\contentsname}{\MakeUppercase\contentsname}

\contentsline {chapter}{Preface}{v}{chapter*.2}%
\contentsline {chapter}{\numberline {1}Overview}{1}{chapter.1}%
\contentsline {part}{I\hspace {1em}\Large The class of Kantorovich Operators}{13}{part.1}%
\contentsline {chapter}{\numberline {2}Dual representations of Kantorovich operators}{15}{chapter.2}%
\contentsline {section}{\numberline {2.1}Skew-linear functionals and entropies}{15}{section.2.1}%
\contentsline {section}{\numberline {2.2}Kantorovich operators associated to skew-linear entropies}{19}{section.2.2}%
\contentsline {section}{\numberline {2.3}A probabilistic reading: conditional large deviations}{23}{section.2.3}%
\contentsline {section}{\numberline {2.4}The convex class of Kantorovich operators}{25}{section.2.4}%
\contentsline {section}{\numberline {2.5}Positively $1$-homogeneous Kantorovich operators}{27}{section.2.5}%
\contentsline {section}{\numberline {2.6}First examples of Kantorovich operators}{31}{section.2.6}%
\contentsline {section}{\numberline {2.7}The non-compact setting}{35}{section.2.7}%
\contentsline {section}{\numberline {2.8}Further comments}{36}{section.2.8}%
\contentsline {chapter}{\numberline {3}Kantorovich operators as functional capacities}{39}{chapter.3}%
\contentsline {section}{\numberline {3.1}Kantorovich operators and Choquet capacities}{39}{section.3.1}%
\contentsline {section}{\numberline {3.2}The Kantorovich envelope of a standard map}{45}{section.3.2}%
\contentsline {section}{\numberline {3.3}The Dellacherie-Kantorovich envelope of a functional capacity}{49}{section.3.3}%
\contentsline {section}{\numberline {3.4}The Choquet-Kantorovich envelope of a functional capacity}{52}{section.3.4}%
\contentsline {section}{\numberline {3.5}Subadditive Kantorovich operators}{54}{section.3.5}%
\contentsline {section}{\numberline {3.6}Further comments}{58}{section.3.6}%
\contentsline {chapter}{\numberline {4}Kantorovich operators on the simplex}{59}{chapter.4}%
\contentsline {section}{\numberline {4.1}A symmetric skew-linear entropy on $X_2$}{59}{section.4.1}%
\contentsline {section}{\numberline {4.2}A homogeneous Kantorovich operator on the simplex}{61}{section.4.2}%
\contentsline {section}{\numberline {4.3}A logarithmic example on $X_n$}{61}{section.4.3}%
\contentsline {section}{\numberline {4.4}Skew-linear entropies associated to convex functions on $\mathbb {R}$}{63}{section.4.4}%
\contentsline {section}{\numberline {4.5}Skew-linear entropies associated to convex functions on $\mathbb {R}^n$}{65}{section.4.5}%
\contentsline {section}{\numberline {4.6}Further comments}{67}{section.4.6}%
\contentsline {chapter}{\numberline {5}Transport and pressure with cost and control}{69}{chapter.5}%
\contentsline {section}{\numberline {5.1}Optimal costs of mass transport as skew-linear entropies}{69}{section.5.1}%
\contentsline {section}{\numberline {5.2}Classical examples of optimal mass transport}{73}{section.5.2}%
\contentsline {section}{\numberline {5.3}Costly pressures as Kantorovich operators}{76}{section.5.3}%
\contentsline {section}{\numberline {5.4}Pressures with cost and control as Kantorovich operators}{83}{section.5.4}%
\contentsline {section}{\numberline {5.5}Further comments}{85}{section.5.5}%
\contentsline {chapter}{\numberline {6}Balayage and homogeneous Kantorovich operators}{87}{chapter.6}%
\contentsline {section}{\numberline {6.1}The classical general theory of balayage}{87}{section.6.1}%
\contentsline {section}{\numberline {6.2}Balayage sets are transitive skew-linear sets}{90}{section.6.2}%
\contentsline {section}{\numberline {6.3}Skew-linear sets are restricted balayage sets}{94}{section.6.3}%
\contentsline {section}{\numberline {6.4}Saturating a skew-linear set into a true balayage set}{97}{section.6.4}%
\contentsline {section}{\numberline {6.5}Balayage with cost}{101}{section.6.5}%
\contentsline {section}{\numberline {6.6}Further comments}{102}{section.6.6}%
\contentsline {part}{II\hspace {1em}\Large Skew-linear and skew-convex entropies}{105}{part.2}%
\contentsline {chapter}{\numberline {7}Skew-linear entropies as optimal cost of balayage}{107}{chapter.7}%
\contentsline {section}{\numberline {7.1}Skew-linear entropies as optimal cost of weak mass transports}{107}{section.7.1}%
\contentsline {section}{\numberline {7.2}Skew-linear entropies as cost of optimal restricted balayage}{113}{section.7.2}%
\contentsline {section}{\numberline {7.3}Subdifferentials of skew-linear entropies}{114}{section.7.3}%
\contentsline {section}{\numberline {7.4}Optimal entropic transport: The Schr\"odinger bridge }{116}{section.7.4}%
\contentsline {section}{\numberline {7.5}Further comments}{116}{section.7.5}%
\contentsline {chapter}{\numberline {8}Operations on skew-linear entropies}{117}{chapter.8}%
\contentsline {section}{\numberline {8.1}The convex cone of skew-linear entropies}{117}{section.8.1}%
\contentsline {section}{\numberline {8.2}Convolution of skew-linear entropies}{121}{section.8.2}%
\contentsline {section}{\numberline {8.3}Entropic regularizations of mass transport }{123}{section.8.3}%
\contentsline {section}{\numberline {8.4}Wasserstein regularization of skew-linear entropies }{127}{section.8.4}%
\contentsline {section}{\numberline {8.5}Further comments}{129}{section.8.5}%
\contentsline {chapter}{\numberline {9}Skew-convex functionals and entropies}{131}{chapter.9}%
\contentsline {section}{\numberline {9.1}Skew-convex functionals and their envelopes}{131}{section.9.1}%
\contentsline {section}{\numberline {9.2}Relative logarithmic entropy is a skew-convex entropy}{138}{section.9.2}%
\contentsline {section}{\numberline {9.3}Convex functions and skew-convex entropies}{143}{section.9.3}%
\contentsline {section}{\numberline {9.4}Skew $\alpha $-convex entropies}{146}{section.9.4}%
\contentsline {section}{\numberline {9.5}Further comments}{149}{section.9.5}%
\contentsline {chapter}{\numberline {10}Inequalities between skew-convex entropies}{151}{chapter.10}%
\contentsline {section}{\numberline {10.1}Bounding skew-convex functionals by entropies}{151}{section.10.1}%
\contentsline {section}{\numberline {10.2}Transport-Entropy inequalities}{155}{section.10.2}%
\contentsline {section}{\numberline {10.3}Maurey-type inequalities}{158}{section.10.3}%
\contentsline {section}{\numberline {10.4}The Prekopa-Leinder theorem and applications}{159}{section.10.4}%
\contentsline {section}{\numberline {10.5}Moment measures}{163}{section.10.5}%
\contentsline {section}{\numberline {10.6}Further comments}{165}{section.10.6}%
\contentsline {part}{III\hspace {1em}Kantorovich operators in Hamiltonian dynamics}{167}{part.3}%
\contentsline {chapter}{\numberline {11}Skew-linear entropies in Lagrangian dynamics}{169}{chapter.11}%
\contentsline {section}{\numberline {11.1}The Lax-Oleinik semi-groups as Kantorovich operators}{169}{section.11.1}%
\contentsline {section}{\numberline {11.2}Eulerian dynamic transport with prescribed end-time}{175}{section.11.2}%
\contentsline {section}{\numberline {11.3}Ballistic skew-linear entropy with fixed end-time}{179}{section.11.3}%
\contentsline {section}{\numberline {11.4}Optimal transport with controlled dynamics and free end times}{183}{section.11.4}%
\contentsline {section}{\numberline {11.5}Eulerian dynamic transport with free end-time}{185}{section.11.5}%
\contentsline {section}{\numberline {11.6}Further comments}{188}{section.11.6}%
\contentsline {chapter}{\numberline {12}Skew-linear entropies in stochastic dynamics}{189}{chapter.12}%
\contentsline {section}{\numberline {12.1}Stochastic mass transports with prescribed end-time}{189}{section.12.1}%
\contentsline {section}{\numberline {12.2}General optimally stopped stochastic transports}{193}{section.12.2}%
\contentsline {section}{\numberline {12.3}Stochastic ballistic transport with fixed end-time}{202}{section.12.3}%
\contentsline {section}{\numberline {12.4}Stochastic mass transport with fixed distribution at all time}{208}{section.12.4}%
\contentsline {section}{\numberline {12.5}Further comments}{209}{section.12.5}%
\contentsline {chapter}{\numberline {13}Dynamic Propagation of Brenier's Mass Transport}{211}{chapter.13}%
\contentsline {section}{\numberline {13.1}The deterministic Bolza duality}{211}{section.13.1}%
\contentsline {section}{\numberline {13.2}Convexity preserving Kantorovich operators}{215}{section.13.2}%
\contentsline {section}{\numberline {13.3}The stochastic Bolza duality}{219}{section.13.3}%
\contentsline {section}{\numberline {13.4}Maximizing the stochastic ballistic cost }{224}{section.13.4}%
\contentsline {section}{\numberline {13.5}Further comments}{227}{section.13.5}%
\contentsline {part}{IV\hspace {1em}\Large Ergodic Properties of Kantorovich Operators}{231}{part.4}%
\contentsline {chapter}{\numberline {14}Mather constants and weak KAM solutions}{233}{chapter.14}%
\contentsline {section}{\numberline {14.1}The Mather constant of a skew-linear entropy}{233}{section.14.1}%
\contentsline {section}{\numberline {14.2}Weak KAM solutions and the Aubry set of a Kantorovich operator }{236}{section.14.2}%
\contentsline {section}{\numberline {14.3}Skew-linear entropies with bounded oscillations}{242}{section.14.3}%
\contentsline {section}{\numberline {14.4}Weak KAM solutions by cone contraction}{244}{section.14.4}%
\contentsline {section}{\numberline {14.5}Skew-linear entropies as large-deviation rates}{246}{section.14.5}%
\contentsline {section}{\numberline {14.6}Further comments}{249}{section.14.6}%
\contentsline {chapter}{\numberline {15}Idempotent Kantorovich operators}{251}{chapter.15}%
\contentsline {section}{\numberline {15.1}Idempotent skew-linear entropies}{251}{section.15.1}%
\contentsline {section}{\numberline {15.2}$\mathcal T$-calibrated functionals}{255}{section.15.2}%
\contentsline {section}{\numberline {15.3}Conjugate idempotent operators}{259}{section.15.3}%
\contentsline {section}{\numberline {15.4}Further comments}{261}{section.15.4}%
\contentsline {chapter}{\numberline {16}Weak KAM operators}{263}{chapter.16}%
\contentsline {section}{\numberline {16.1}Weak KAM operators and Measure-level Aubry sets}{263}{section.16.1}%
\contentsline {section}{\numberline {16.2}Weak KAM operators when $T$ and $T^2$ are comparable}{267}{section.16.2}%
\contentsline {section}{\numberline {16.3}Weak KAM operators for oscillation contractions}{269}{section.16.3}%
\contentsline {section}{\numberline {16.4}Balanced skew-linear entropies}{270}{section.16.4}%
\contentsline {section}{\numberline {16.5}Power bounded Kantorovich contractions}{275}{section.16.5}%
\contentsline {section}{\numberline {16.6}Mather theory for continuous skew-linear entropies}{279}{section.16.6}%
\contentsline {section}{\numberline {16.7}Amenable skew-linear entropies}{284}{section.16.7}%
\contentsline {section}{\numberline {16.8}Weak KAM theory for positively $1$-homogeneous operators}{288}{section.16.8}%
\contentsline {section}{\numberline {16.9}Further comments}{290}{section.16.9}%
\contentsline {part}{V\hspace {1em}\Large Weak KAM operators in classical analysis}{291}{part.5}%
\contentsline {chapter}{\numberline {17}Weak KAM operators in linear and convex analysis}{293}{chapter.17}%
\contentsline {section}{\numberline {17.1}Weak KAM operators associated to Markov operators}{293}{section.17.1}%
\contentsline {section}{\numberline {17.2}Concavification weak KAM operators on open domains}{295}{section.17.2}%
\contentsline {section}{\numberline {17.3}Weak KAM concavification operators on general domains}{298}{section.17.3}%
\contentsline {section}{\numberline {17.4}Choquet theory as a Mather theory}{302}{section.17.4}%
\contentsline {section}{\numberline {17.5}Further comments}{303}{section.17.5}%
\contentsline {chapter}{\numberline {18}Weak KAM operators in complex analysis}{305}{chapter.18}%
\contentsline {section}{\numberline {18.1}Plurisuperharmonic envelopes on open domains}{305}{section.18.1}%
\contentsline {section}{\numberline {18.2}Jensen barycenters and $PSH$-convex sets}{309}{section.18.2}%
\contentsline {section}{\numberline {18.3}Pluri-superharmonic envelopes as weak KAM operators}{312}{section.18.3}%
\contentsline {section}{\numberline {18.4}Holomorphic images of complex Brownian motion}{313}{section.18.4}%
\contentsline {section}{\numberline {18.5}Further comments}{317}{section.18.5}%
\contentsline {chapter}{\numberline {19}Weak KAM operators and optimal Brownian stopping}{319}{chapter.19}%
\contentsline {section}{\numberline {19.1}The cone of superharmonic functions}{319}{section.19.1}%
\contentsline {section}{\numberline {19.2}Superharmonic envelopes as weak KAM operators}{323}{section.19.2}%
\contentsline {section}{\numberline {19.3}Skorokhod embeddings in Brownian motion }{327}{section.19.3}%
\contentsline {section}{\numberline {19.4}Optimal Brownian Transport with cost}{329}{section.19.4}%
\contentsline {section}{\numberline {19.5}Further comments}{333}{section.19.5}%
\contentsline {part}{VI\hspace {1em}\Large Mather-Aubry Theory and its variations}{335}{part.6}%
\contentsline {chapter}{\numberline {20}Mather--Aubry theory for free-energy operators}{337}{chapter.20}%
\contentsline {section}{\numberline {20.1}The Aubry set of a free-energy operator}{337}{section.20.1}%
\contentsline {section}{\numberline {20.2}Existence of continuous weak KAM solutions}{339}{section.20.2}%
\contentsline {section}{\numberline {20.3}The Schr\"odinger semigroup}{342}{section.20.3}%
\contentsline {section}{\numberline {20.4}The Sinkhorn operator}{343}{section.20.4}%
\contentsline {section}{\numberline {20.5}The Ruelle free-energy operator}{344}{section.20.5}%
\contentsline {section}{\numberline {20.6}Further comments}{346}{section.20.6}%
\contentsline {chapter}{\numberline {21}Mather theory for general optimal mass transports}{347}{chapter.21}%
\contentsline {section}{\numberline {21.1}Idempotent optimal mass transport}{347}{section.21.1}%
\contentsline {section}{\numberline {21.2}Discrete Mather theory for mass transport}{349}{section.21.2}%
\contentsline {section}{\numberline {21.3}Mather theory for a continuous semi-group of transports}{353}{section.21.3}%
\contentsline {section}{\numberline {21.4}The Sinkhorn operator and the Schr\"odinger bridge}{354}{section.21.4}%
\contentsline {section}{\numberline {21.5}A mixing kernel with cost}{356}{section.21.5}%
\contentsline {section}{\numberline {21.6}Further comments}{357}{section.21.6}%
\contentsline {chapter}{\numberline {22}Mather Theory in Hamiltonian dynamics}{359}{chapter.22}%
\contentsline {section}{\numberline {22.1}Role of the stationary Hamilton-Jacobi equation}{359}{section.22.1}%
\contentsline {section}{\numberline {22.2}Viscosity solutions of Hamilton--Jacobi equations}{361}{section.22.2}%
\contentsline {section}{\numberline {22.3}Weak KAM solutions and the Fathi approach}{364}{section.22.3}%
\contentsline {section}{\numberline {22.4}The Bernard-Buffoni approach}{366}{section.22.4}%
\contentsline {section}{\numberline {22.5}A brief description of the Aubry-Mather theory}{370}{section.22.5}%
\contentsline {section}{\numberline {22.6}Stochastic weak KAM on the Torus}{372}{section.22.6}%
\contentsline {section}{\numberline {22.7}Further comments}{375}{section.22.7}%
\contentsline {chapter}{\numberline {23}Mather theory in Ergodic Optimization}{377}{chapter.23}%
\contentsline {section}{\numberline {23.1}The standard model}{377}{section.23.1}%
\contentsline {section}{\numberline {23.2}The general holonomic model}{382}{section.23.2}%
\contentsline {section}{\numberline {23.3}Ergodic optimization in symbolic dynamics}{387}{section.23.3}%
\contentsline {section}{\numberline {23.4}A stochastic holonomic setting}{388}{section.23.4}%
\contentsline {section}{\numberline {23.5}Further comments}{392}{section.23.5}%
\contentsline {chapter}{\numberline {24}Controlled and risk-sensitive weak KAM theory}{395}{chapter.24}%
\contentsline {section}{\numberline {24.1}The cost-and-control operator and its examples}{395}{section.24.1}%
\contentsline {section}{\numberline {24.2}Existence of weak KAM solutions}{397}{section.24.2}%
\contentsline {section}{\numberline {24.3}The large-deviation reading}{400}{section.24.3}%
\contentsline {section}{\numberline {24.4}Further comments}{400}{section.24.4}%
\contentsline {chapter}{\numberline {25}Appendix: Towards a notion of skew-multilinear entropies}{403}{chapter.25}%
\contentsline {section}{\numberline {25.1}The $\mathcal T$-Legendre transform and skew-multilinear entropies}{403}{section.25.1}%
\contentsline {section}{\numberline {25.2}Skew-multilinear entropies as weak transport}{407}{section.25.2}%
\contentsline {section}{\numberline {25.3}Operations on skew-multilinear entropies}{410}{section.25.3}%
\contentsline {section}{\numberline {25.4}On an alternate definition}{413}{section.25.4}%
\contentsline {chapter}{\numberline {26}References}{415}{chapter.26}%

\makeatother

\end{document}